\documentclass[11pt,draftcls,onecolumn]{IEEEtran}
\usepackage{amsfonts}
\usepackage{amsmath}
\usepackage{amssymb}
\usepackage{algorithm}
\usepackage{booktabs}
\usepackage{graphicx}
\usepackage{color}
\usepackage[dvipsnames]{xcolor}
\usepackage{enumerate}
\usepackage{amsbsy}
\usepackage{booktabs}
\usepackage{subcaption}
\usepackage{url}
\usepackage{cite}
\usepackage{cuted,tcolorbox,lipsum,tikz,pstricks} 
\usepackage{wrapfig,color}
\usepackage{pst-node,graphicx}
\SpecialCoor
\definecolor{carmine}{rgb}{0.59, 0.0, 0.09}

\newtheorem{remark}{Remark}

\newcommand{\real}{\text{Re}}
\newcommand{\imag}{\text{Im}}
\renewcommand{\i}{{\text{i}}}

\newcommand{\ped}[1]{{_{\mathrm{#1}}}}

\newcommand{\uu}{\mathbf{u}}
\newcommand{\xx}{\mathbf{x}}

\newcommand{\WW}{\mathbf{W}}
\newcommand{\dd}{\boldsymbol{\delta} }
\renewcommand{\aa}{\boldsymbol{\alpha}}
\newcommand{\DD}{\boldsymbol{\Delta} }

\definecolor{green}{rgb}{0,0.7,0}

\IEEEoverridecommandlockouts

\title{Algorithms for Optimal AC Power Flow in the Presence of Renewable Sources}

\begin{document}

\author{\IEEEauthorblockN{Mohammadreza Chamanbaz\IEEEauthorrefmark{1},
Fabrizio Dabbene\IEEEauthorrefmark{2}, and
Constantino Lagoa\IEEEauthorrefmark{3}}\\
\IEEEauthorblockA{\IEEEauthorrefmark{1}Arak University of Technology, Iran\\
\IEEEauthorrefmark{2}CNR-IEIIT, Politecnico di Torino\\
\IEEEauthorrefmark{3}The Pennsylvania State University, USA}
\thanks{This work was supported in part by the CNR International Joint Lab COOPS, the National Science Foundation under grant \mbox{CNS-1329422}, the Singapore National Research Foundation (NRF) grant
under the ASPIRE project, grant No NCR-NCR001-040 and Iran's National Elites Foundation.}}


\maketitle
	
\begin{abstract}
This chapter presents recent solutions to the optimal power flow (OPF) problem in the presence of renewable energy sources (RES), {such} as solar photo-voltaic and wind generation.
After introducing the original formulation of the problem, arising from the combination of economic dispatch and power flow, we provide a brief overview of the different solution methods proposed in the literature to solve it. Then, we explain the main difficulties arising from the increasing RES penetration, and the ensuing necessity of deriving robust solutions.
Finally, we present the state-of-the-art techniques, with a special focus on recent methods we developed, based on the application on randomization-based methodologies.
\end{abstract}

\section{Introduction}


The main purpose of an electrical network is to transfer electrical power from the generation sources to the consumers. In order to reliably and effectively perform this task,   efficient ways of doing power  transfer need to be developed that i) satisfy demand,  ii) minimize the cost of power generation and iii)~do not violate safety constraints. This is the main objective of the so-called Optimal Power Flow (OPF) problem; i.e., to determine how much power should each generator in the network produce so that the safety limits of both generators and lines are not exceeded, demand is satisfied  and power is generated at the lowest possible cost.


The OPF problem, whose formulation can be traced back to the sixties  \cite{Carpentier:1962aa}, is  relatively easy to state but quite difficult to solve. More precisely, once the equations that describe  power balance and  safety constraints are obtained, one realizes that the problem of determining the optimal power generation policy is a very complex, non-convex one for which off-the-shelf solvers very often do not provide satisfactory solutions.


The OPF problem becomes even more complex if renewable energy sources (RESs) are present in the network. 
Modern power grids are characterized by increasing penetration of RES, such as solar photovoltaic and wind power. 
This trend is expected to increase in the near future, as also testified by strict commitments to large renewable power penetration being made by major countries worldwide; e.g., see  \cite{US-DOE2008,EWEA2009,SET2016,GWEC2016}. 


While the advantages of renewable energy in terms of environmental safeguard are indisputable, its introduction does not come without a cost. Indeed, renewable energy generation technologies are highly variable and not fully dispatchable, thus imposing novel challenges to the existing power system operational paradigm. As discussed in e.g. \cite{Bienstock2014}, when  uncontrollable resources fluctuate, classical optimal power flow  solutions can provide very inefficient power generation policies, that result in line overloads and, potentially,  cascading outrages.

Despite the increasingly larger investments, which are costly and subject to several regulatory and policy limitations, the frequency and scale of power outages are steadily growing. This situation clearly shows that a strategy only based on investments in technological improvements of the transmission lines and controllable generation capacity---as those discussed e.g.\ in \cite{borkowska1974,Conejo2010}---is not 
sufficient anymore. Instead, radically new dispatch philosophies need to be devised, able to cope with the increasing \textit{uncertainty}, due to unpredictable fluctuations in renewable output and time-varying loads. 

Indeed, classical OPF dispatch---employed to design generator output setpoints in order to meet demand at minimum cost, without violating operating limitations of generators and transmission lines---is typically computed based on simple predictions of expected loads and generation levels for the upcoming time window. Although these predictions can be fairly precise for the case of traditional generators and loads---thus allowing OPF to achieve remarkable reliability of operation in the presence of normal fluctuations---they may be highly unreliable in the case of renewable generators, thus explaining its failure in these latter situations. As a consequence, recent regulatory initiatives, such as Federal Energy Regulatory Commission (FERC) Orders 764 \cite{FERC2012} and 890 \cite{FERC2007} have identified the need for a new generation of operating protocols and decision-making tools for the successful integration of renewable generation.

It follows that one of the major problems in today's power grids is the following \cite{OurTCNS2018}: \textit{Given the high level of uncertainty introduced  by renewable energy sources, design a dispatch policy that {\it i)} minimizes generation costs and {\it ii)} has a very small risk of violating generation and transmission constraints. In other words, one would like to design an optimal dispatch policy with very low risk of network failure.}


Although a lot of attention has been devoted to the classical OPF problem (without uncertainty), only a limited number of published articles address---in a systematic way---the case where uncertain power generation is present like the one described above.

%
%
%

\section{Nominal Optimal Power Flow}

\subsection{Formulation of nominal OPF}

The main purpose of a power network is to transfer electrical power from generators to consumers (loads). Optimal power flow  manages the network and controls all controllable parameters such as active power generation, transformer taps, capacitor banks, etc. Almost all OPF formulations are derived from the classical formulation of Carpentier in 1962, see  \cite{Carpentier:1962aa}. The most common form of OPF problem seeks a solution minimizing the total electrical generation cost while at the same time ensuring that the power network remains in its safe operating region. In fact, safe operation is enforced by adding enough constraints, which form the feasible set of the mathematical optimization problem. The optimal point is then selected as the point having the smallest cost in the feasible set. 

Power network is usually modeled as a set of nodes $\mathcal{N} \doteq \{1,2,\ldots,n\}$ and a set of edges $\mathcal{L}\subset\mathcal{N}\times \mathcal{N}$. Each node represents a bus and each edge represents an electrical line connecting two buses in the power network. Therefore, we use bus and node interchangeably, and line and edge interchangeably. If $(i,j) \in \mathcal{L}$, this means that there is a line connecting buses $i$ and $j$.  Each bus $k$ in the network is characterized by two complex variables: complex voltage $V_k = |V_k|\angle \theta_k$---where $|V_k|$  and $\angle \theta_k$ are  the magnitude and angle of the complex voltage $V_k$ respectively---and complex power $S_k = P_k+Q_k\i$---where $P_k$ and $Q_k$ are the active and reactive power in bus $k$. 

Based on this description, buses in the power network can be categorized into three groups.
The first bus type is the slack bus, usually denoted as bus $0$. In most cases, the voltage magnitude and phase are fixed at the slack bus---typically $V_0 = 1\angle0$---whereas active and reactive generator powers are variables. Slack bus is considered as \emph{reference} bus; its role is to balance the active and reactive power in the power grid. The slack bus must include a generator. 

The second bus type is the generator bus. In a generator bus active power $P^G_k$ and voltage magnitude $|V_k|$ are given, while reactive power $Q^G_k$ and voltage angle $\theta_k$ are variable (we assume for simplicity of notation that no more than one generator is present on each generator bus). A generator bus is also called PV bus. We denote by $\mathcal{G}$ the set of generator buses. The active and reactive powers of a generator are usually restricted by some operational limits:
\begin{subequations}
\makeatletter
        \def\@currentlabel{PL}
        \makeatother
        \label{eq:power constraints}
        \renewcommand{\theequation}{PL.\arabic{equation}}
\begin{align}
\label{eq:power constraints active generation bound}
& P_{k\,\min}\leq P^G_{k}\leq P_{k\,\max},\quad \forall k \in \mathcal{G} \\
\label{eq:power constraints reactive generation bound}
& Q_{k\,\min}\leq Q^G_{k}\leq Q_{k\,\max},\quad \forall k \in \mathcal{G}
\end{align}
\end{subequations} 
for given limits $P_{k\, \min},\,P_{k\, \max}$ and $Q_{k\, \min},\,Q_{k\, \max}$. 
We refer to \eqref{eq:power constraints} as \textit{(Active and Reactive) Power Limits}.

The third type is the load bus---also called PQ bus---in which the complex power $P^L_k+Q^L_k\i$ is specified and the complex voltage $|V_k|\angle \theta_k$ is variable.  We denote by $\mathcal{D}$ the set of load buses.

We adopt the convention of choosing bus $0$ as the slack bus and the first $n_g$ buses  as generator buses, i.e.\ $\mathcal{G}=\{1,2,\ldots,n_g\}$. The remaining buses $\mathcal{D}=\{n_g+1,\ldots,n\}$ are hence load buses.

Clearly, it is of paramount importance that the generation levels are allocated to the generating units $j\in\mathcal{G}$, so that the system load is supplied entirely and most economically, while guaranteeing the Power Limits. This problem is usually referred to (static) Economic Dispatch (ED) \cite{book-PF4}, and represents the simplest optimization problem one can consider in designing the power grid operations. The objective of ED is to calculate, for a single period of time, the output active power of every generating unit so that all the demands are satisfied at the minimum cost, while satisfying different technical constraints of the network and the generators. To this end, to each generating unit a unique production cost is associated, in the form of a cost function  $f_k,\,k\in\mathcal{G}$. Then, in its simplest form, the problem writes as follows

\noindent{\bf Economic Dispatch Problem}
\begin{eqnarray} \label{eq:ED}
&\underset{P^G_1\,  \cdots \,P^G_{n_g}
}{\text{minimize}}
&\sum_{k\in\mathcal{G}} f_k(P^G_{k})\\ \nonumber
&\text{subject to }
&\sum_{k\in\mathcal{G}} P^G_{k}=P^D\nonumber\\
&&P_{k\,\min}\leq P^G_{k}\leq P_{k\,\max},\quad \forall k \in \mathcal{G},\nonumber
\end{eqnarray}
where $P^D$ represents the total demand. The equality constraint in \eqref{eq:ED} represents the demand--supply balance constraint, while the inequality represents the Active Power Limits \eqref{eq:power constraints active generation bound}. Hence, the ED problem is to minimize cost while guaranteeing that supply meets demand. 

It should be remarked that, while ED accounts for generator limits, it \textit{does not consider network} constraints. Indeed, the power network needs to satisfy the so called \textit{Power Flow equations}, which  are in fact nothing other than Kirchhoff's circuit laws applied to the power network. In general, there are two methods to write the power flow equations: i) bus injection model, which is the most compact form and ii) branch flow model. In the present paper, we use branch flow model,
since it presents several advantages in terms of convex relaxation, see Section~\ref{sec:relax}. It is worth highlighting that the two models are equivalent. 
To define the Power Flow equations, consider the line $(\ell,m)\in \mathcal{L}$, which is the line connecting buses $\ell$ and $m$.  Let  
$Y_{\ell m}=G_{\ell m}+B_{\ell m}\i$ be the (complex) admittance of the line and $V_k=|V_k|\angle{\theta_k}$ be the (complex) voltage at bus $k$. Then, the following balance equations should be satisfied at all times

\noindent{\bf Power Flow Problem}
\begin{subequations}
\makeatletter
        \def\@currentlabel{PF}
        \makeatother
        \label{eq:balance equations}
        \renewcommand{\theequation}{PF.\arabic{equation}}
\begin{align}
\label{eq:active balance equation}
& P^G_{k} -  P^L_{k} = \sum_{\ell\in\mathcal{N}_k} \real\left\{V_k(V_k-V_l)^*Y_{k\ell}^*\right\}, \quad \forall k \in \mathcal{N}\\
\label{eq:OPF reactive balance equation}
& Q^G_{k} -  Q^L_{k} = \sum_{\ell\in\mathcal{N}_k} \imag\left\{V_k(V_k-V_l)^*Y_{k\ell}^*\right\},\quad \forall k \in \mathcal{N},
\end{align}
\end{subequations}
where $\mathcal{N}_k$ is the set of all neighboring buses directly connected to bus $k$, and  $P^L_{k}$ and $Q^L_{k}$ denote the (known) active and reactive loads of the network at bus $k$\footnote{{Note that, by convention,  $P^G_{k}$ and $Q^G_{k}$ are set to zero in non-generator nodes.}}. Equations \eqref{eq:balance equations} simply mean that the active and reactive power at bus $k$ need to be balanced. The Power Flow Problem hence amounts to finding a feasible solution to equations \eqref{eq:balance equations}.

Finally, in general, the network should satisfy requirements also in terms of bus voltages. In particular, the following \textit{Voltage Limits} are usually imposed:
\begin{subequations}
\makeatletter
        \def\@currentlabel{VL}
        \makeatother
        \label{eq:voltage constraints}
        \renewcommand{\theequation}{VL.\arabic{equation}}
\begin{align}
\label{eq:voltage bound 1}
&V_{k\,\min}\leq |V_k|\leq V_{k\,\max}, \quad \forall k \in \mathcal{N}, \\
\label{eq:voltage bound 2}
& |V_\ell-V_m|\le \Delta V_{\ell,m}^{\max},
 \quad \forall (\ell,m) \in \mathcal{L}.
\end{align}
\end{subequations}

Inequality \eqref{eq:voltage bound 1} restricts the magnitude of the voltage at bus $k$ and \eqref{eq:voltage bound 2} restricts the amount of power carried by the line $(l,m)\in\mathcal{L}$. Each line in the power network has a limitation in terms of the amount of complex power passing through the line. The constraint \eqref{eq:voltage bound 2} represents one of the most important constraints in the power network: frequent violation of \eqref{eq:voltage bound 2} can lead to overheating of the line and consequently line tripping.
In \cite{Lavaei:2012}, the constraints in \eqref{eq:voltage bound 2} have been proven to be practically equivalent to the more classical bound
\[
 \left | V_l(V_l-V_m)^*y_{lm}^*\right | \le S_{lm\, \max},\quad\forall(l,m)\in\mathcal{L},
\]
where $S_{lm\, \max}$ is the maximum apparent power flow allowed to pass through the line $(l,m)\in\mathcal{L}$.

The  OPF problem arises from combining Economic Dispatch and Power Flow problems:  the goal is to minimize generation cost subject to the Power Flow constraints \eqref{eq:balance equations} and to the operational constraints \eqref{eq:power constraints} and \eqref{eq:voltage constraints}. 
Formally, the \textit{AC-OPF Problem} can be stated as\\
\noindent{\bf AC - Optimal Power Flow}
\begin{subequations} 
        \makeatletter
        \def\@currentlabel{AC-OPF}
        \makeatother
        \label{eq:OPF Original}
        \renewcommand{\theequation}{AC-OPF}
\begin{align}
\underset{P^G_1\,  \cdots \,P^G_{n_g},
Q^G_1\,  \cdots \,Q^G_{n_g},
V_1\,  \cdots \,V_{n}
}{\text{minimize}} \ \ 
&\sum_{k\in\mathcal{G}} f_k(P^G_{k})\\ \nonumber
\text{subject to:} \ \  
&\text{\it Power Flow  Equations }\eqref{eq:balance equations},\nonumber\\
&\text{\it Power Limits }\eqref{eq:power constraints},\nonumber\\
&\text{\it Voltage Limits }\eqref{eq:voltage constraints}.\nonumber
\end{align}
\end{subequations}
The individual generation cost  $f_k$, representing the cost of generation for generator $k$, is usually chosen as a quadratic function. 
The formulation of the nominal AC-OPF problem in \eqref{eq:OPF Original} is very classical and can be traced back to the seminal paper of Carpentier \cite{Carpentier:1962aa}. We remark that problem  \eqref{eq:OPF Original} is a difficult optimization problem. Specifically, it is  non-convex
 due to the nonlinear equations~\eqref{eq:balance equations}, which contain quadratic equalities involving the voltages. The nonconvexity of the OPF problem has sparked several research directions for the practical solution of the problem. This is discussed in the next section.\\
 

\begin{tcolorbox}[title= \sf \scriptsize \textbf{Small Note on Convexity},colframe=carmine!10,
colback=carmine!10,
coltitle=black,  
]
\sf\scriptsize

The general form of an optimization problem is 
\begin{eqnarray*} \label{eq:general convex problem}
&\underset{x
}{\text{minimize}}
&c(x)\\ \nonumber
&\text{subject to }
&g_i(x)\leq 0, i = 1,\ldots,m, \nonumber
\end{eqnarray*}
where $x \in \mathbb{R}^n$ is the vector of optimization variables, $c(x)$ is the scalar \textit{objective function} or \textit{cost function} and $g_i(x)\leq0,i = 1,\ldots,m$ are scalar constraints  defining the feasible set of $x$. 
Convexity plays a crucial role in problem \eqref{eq:general convex problem}. Simply speaking, a function is convex if it curves up. More precisely, a function $f:\mathbb{R}^n\rightarrow \mathbb{R}$ is convex if its domain is a convex set and for any $x,y$ in the domain of $f$, the following inequality holds
\[
f(\theta x+(1-\theta)y)\leq \theta f(x) + (1-\theta) f(y), \qquad0\leq\theta\leq 1.
\]
The main property of a convex function is that its first order Taylor approximation is a global underestimator for the function $f$. Optimization problem above is convex if the cost function $c(x)$ and all constraints $g_i(x),i=i,\ldots,m$ are convex in the decision variable $x$. There is a wealth of numerical algorithms available to solve a convex problem. On the contrary, if a problem is non-convex, most algorithms are only able to find a local minimum and the global minimum is usually very difficult to find. 
\end{tcolorbox}


\begin{remark}[Unit commitment]
A problem that is closely related to the  problem \eqref{eq:OPF Original} described above is the so-called \textit{Unit Commitment Problem}. This is an extension of OPF where i) a time interval is considered, not just a specific time instant and ii) costs of generator operation such as start-up and shutdown costs are included. To keep the presentation simple and accessible, we do not discuss this problem in this chapter, and refer the interested reader to, e.g.,~\cite{Saravanan2013,LUJANOROJAS2016}. We point out that the different approaches to OPF discussed in this chapter can be extended to deal with this more complex problem.
\end{remark}

\section{Solution Approaches and Relaxations for OPF}

The problem of dealing with the nonlinear equality and inequality constraints in~\eqref{eq:balance equations} may be addressed in two main ways. On the one end, several global optimization approaches have been developed to directly solve the nonconvex problem. Clearly, these approaches suffer from the typical drawbacks of global optimization tools: they are guaranteed to converge only to a local minimum and, in general, they may show a rather high computational complexity, making them hardly scalable to large networks. We briefly review these approaches in Section \ref{sec:opttools}.

On the other end, different simplifying assumptions have been introduced in the literature aimed at deriving  computable approximations of the nonconvex problem.
These techniques are briefly reviewed in Section \ref{sec:DC} (devoted to the so-called DC-approximations) and Section \ref{sec:relax}, which concentrates on a more modern approach based on semidefinite programming (SDP) relaxations.

\subsection{Optimization Methods}
\label{sec:opttools}
Optimization methods to solve the OPF problem \eqref{eq:OPF Original}  can be divided into two different groups: i) deterministic methods and ii) heuristic search techniques. 

\subsubsection{Deterministic Optimization Methods}

Almost all deterministic techniques in optimization literature have been applied to OPF problem. Some of the deterministic methods for solving OPF problem are stated next.

\begin{enumerate}
	\item 
	\emph{Gradient Methods:}
	Gradient methods were 
	the first techniques used in $1960$s to solve the OPF problem. These techniques can be classified into three main subclasses: i) Reduced gradient methods which were first introduced in \cite{wolfe1967methods} in 1967, and applied for the first time to OPF problem in \cite{dommel1968optimal}, ii) Conjugate gradient methods \cite{burchett1982developments} which is an extension of the reduced gradient method, and iii) Generalized reduced gradient methods, which are a further generalization of the reduced gradient methods, see \cite{abadie1969generalization} and \cite{peschon1972optimal} for a detailed explanation of this approach. The main feature of gradient based methods is that they use first order derivative to sequentially move towards the local minimum point leading to slow progress compared to methods using higher order derivatives, e.g. Newton, quasi Newton, etc. 
	
	\item
	\emph{Newton's Method:} 
	The second class of algorithms used for solving OPF problem is Newton's method, which uses the second order derivative (Hessian) in order to improve convergence speed \cite{deuflhard2011newton}.  Newton methods require Hessian matrix and its inverse at each iteration of the algorithm, which is computationally demanding. Nevertheless, in most OPF problems the Hessian matrix can be computed analytically and efficient numerical methods can be employed to reduce the computational burden associated with matrix inversion. Quasi-Newton methods are a class of algorithms which rather approximate the Hessian matrix and perform faster than Newton methods in some cases \cite{housos1982sparse}. 
	\item
	\emph{Interior Point Methods:}
	The third class of optimization algorithms used in OPF are Interior Point Methods (IPMs). This class of methods limits the search space to the interior of the feasible region by augmenting the objective with some barrier terms \cite{karmarkar1984new}.  IPMs need few iterations for convergence, although each iteration is computationally demanding. They are the most widely used algorithms for solving OPF problem. There are a number of different extensions to the IPM, such as primal-dual interior point methods \cite{granville1994optimal}, predictor-corrector primal-dual interior point methods \cite{torres1998interior} and trust region techniques \cite{min2005trust}. 
	
	\item
	\emph{Other Methods:} 
	Other algorithms such as simplex methods \cite{Dantzigsimplex}, sequential linear programming \cite{griffith1961nonlinear} and sequential quadratic programming \cite{bazaraa2013nonlinear} are also used for solving different versions of OPF problem.   
	
\end{enumerate}

\subsubsection{Heuristic Optimization Methods}
As mentioned earlier, the OPF optimization problem is non-convex. Heuristic methods perform better in the face of a non-convex problem as they are designed to avoid local minima. Most heuristic methods borrow ideas from natural phenomena such as biological evolution, school of fish and flock of birds. There are a number of heuristic---also known as random search or stochastic---methods used for solving OPF problem. Some of them are discussed next.

\begin{enumerate}
	\item 
	\emph{Evolutionary Algorithms:} These class of algorithms mimic biological evolution. Candidate solutions are individual population which are progressively improved in terms of cost function by applying various biological rules such as selection, recombination, and mutation. A wide class of algorithms fall within the  category of evolutionary algorithms. Genetic algorithms \cite{Goldberg:1989:GAS:534133}, evolutionary programming \cite{fogel2006evolutionary}, artificial immune systems \cite{kephart1994biologically}  and differential evolution \cite{storn1997differential} are a number of evolutionary algorithms used for solving different versions of the OPF problem; see e.g. \cite{Abido_environmental,ABOUELELA2010878} for an implementation of evolutionary algorithms on economic power dispatch and OPF respectively. 
	
	\item
	\emph{Artificial Neural Network:} 
	This methods works based on similar principles as biological neural network constituting human/animal brain \cite{dreyfus2005neural}. It is parallel in nature and has the ability to learn as the new data is injected to the algorithm. 
	
	\item
	\emph{Particle Swarm Optimization:} 
	This method solves an optimization problem by generating a large number of candidate solutions (a.k.a population or particles). Each candidate solution is then updated using a formula incorporating the position and velocity of the  particles \cite{kennedy2011particle,abido2002optimal}. This method borrows ideas from natural swarms such as flocks of birds and schools of fish.  
	
	\item
	\emph{Simulated Annealing:}
	At each iteration, the algorithm selects a candidate solution close to the current one and evaluates its ``goodness'' in terms of the cost function. There are two probabilities evolving during the execution of the algorithm: i) probability of moving to a worse solution which is progressively moved towards zero  and ii) probability of moving to a better solution which is either kept to $1$ or changed to a positive value \cite{kirkpatrick1983optimization}. Selecting a worse candidate solution prevents the algorithm from getting stuck in a local minimum.  
	
	\item
	\emph{Other Methods:}
	Other heuristic optimization algorithms such as ant colony optimization \cite{maniezzo1992distributed}, bacterial foraging \cite{passino2002biomimicry}, Tabu search \cite{Abido_Tabu_Search}, and chaos optimization algorithm \cite{jiang1998optimizing} have been applied to OPF problem.
	
\end{enumerate}

\subsection{DC OPF}
\label{sec:DC}
The most common approximation approach to deal with the AC Power Flow equations consists in linearizing them by introducing a series of simplifying restrictions regarding voltage magnitudes, voltage angles, admittances, and reactive power, justified by operational considerations under normal operating conditions. Many different versions of the so-called DC approximation model exist, see for instance the books \cite{book-PF1,book-PF2,book-PF3,book-PF4}. Clearly, a comprehensive review of all the different variants falls outside the scope of this work: the reader is referred to the recent works \cite{Stott:2009,Coffrin:2014} which present  an in-depth discussion together with a modern re-interpretation of the problem. 
For brevity, we here review the simplest and most popular variant of the DC approximation.
To this end, remembering that $V_k = |V_k|\angle \theta_k$ and $Y_{\ell m}=G_{\ell m}+B_{\ell m}\i$, we first rewrite equations \eqref{eq:balance equations} as follows
\begin{align*}
& P^G_{k} -  P^L_{k} = \sum_{\ell\in\mathcal{N}_k} |V_k||V_\ell|\left[G_{k\ell}\cos(\theta_k-\theta_\ell)+B_{k\ell}\sin(\theta_k-\theta_\ell)\right], \quad \forall k \in \mathcal{N}\\
& Q^G_{k} -  Q^L_{k} = \sum_{\ell\in\mathcal{N}_k} |V_k||V_\ell|\left[G_{k\ell}\cos(\theta_k-\theta_\ell)-B_{k\ell}\sin(\theta_k-\theta_\ell)\right],
\quad \forall k \in \mathcal{N}.
\end{align*}
In the classical DC approximation, the following assumptions are made:
\begin{enumerate}[i)]
    \item the susceptance is large relative to the conductance, i.e. $|G_k| << |B_k|$;
    \item the phase angle difference is small enough to ensure $\sin(\theta_k-\theta_\ell)\approx (\theta_k-\theta_\ell)$ and
     $\cos(\theta_k-\theta_\ell)\approx 1$ (the term disappears because $G_{kl}$ small);
    \item the voltage magnitudes $V_k$ are close to 1.0 and do not vary significantly. 
\end{enumerate}
Under these assumptions, equations \eqref{eq:balance equations} reduce to
\begin{equation}
P^G_{k} -  P^L_{k} = \sum_{\ell\in\mathcal{N}_k} B_{k\ell}(\theta_k-\theta_\ell), \quad \forall k \in \mathcal{N}.
\end{equation}
Although the DC approximation leads to a tractable problem, the solution obtained using this method is in general sub-optimal and, more importantly,  it may not be feasible, in the sense that it may not satisfy the original nonlinear power flow equations. In this latter situation, the network operator has to rerun the optimization modifying the constraints adopting some ad-hoc heuristics, with in general no  \textit{a-priori} guarantee of convergence. This renders this solution not suitable for large systems. 
Also, as noted in \cite{Coffrin:2014},  the fact that DC approximation fixes voltage magnitudes and ignores reactive power, makes the solution not applicable in several important practical situations.

\subsection{Convex relaxation of OPF}
\label{sec:relax}
In recent years, a more sophisticated relaxation has been introduced in \cite{Bai:08,Lavaei:2012} for circumventing the non-convexity associated with the optimal power flow problem \eqref{eq:OPF Original}. This relaxation stems from the consideration that the source of nonconvexity is due to non-linear {(quadratic)} terms $V_kV_l$'s appearing in the equality constraint \eqref{eq:balance equations} and in the inequality constraint \eqref{eq:voltage constraints}. It can be noted that  the quadratic constraints can be reformulated as linear ones by introducing a new variable 
\[
\mathbf{W}=\mathbf{V}\mathbf{V}^*
\]
where $\mathbf{V} \doteq[V_1,\ldots,V_n]^T$ is the vector of complex bus  voltages. In order to replace $\mathbf{V}\mathbf{V}^*$ with the new variable $\mathbf{W}$, two additional constraints need to be included: 
\begin{enumerate}[i)]
    \item the matrix $\mathbf{W}$ needs to be positive semi definite 
i.e. the following \textit{Positivity Constraint} should hold
\begin{equation}
\label{positivity_constr} 
\mathbf{W}\succeq 0,
\end{equation}
\item the rank 
 of $\mathbf{W}$ should be one, i.e.\ the following \textit{Rank Constraint} should hold
\begin{equation}
\label{positity_constr} 
\text{rank}\{\mathbf{W}\}=1. 
\end{equation}
\end{enumerate}

\begin{tcolorbox}[title= \sf \scriptsize \textbf{Positive Semidefinite Matrix},colframe=carmine!10,
colback=carmine!10,
coltitle=black,  
]
\sf\scriptsize
A conjugate symmetric matrix $M=M^* \in C^{n \times n}$ is said to be a positive semidefinite matrix, denoted by
$M \succeq 0
$ if \textit{for all} $x\in C^n$ it holds
\[
x^* M x \geq 0.
\]
We should note at this point that the set of all positive definite matrices 
$\mathcal{S}^{++}_0 = \{ M=M^* \in C^{n \times n} : M \succeq 0\}
$ is a convex set.

\end{tcolorbox}

\begin{tcolorbox}[title= \sf \scriptsize \textbf{Rank of a Matrix},colframe=carmine!10,
colback=carmine!10,
coltitle=black,  
]
\sf\scriptsize

Given a matrix $M \in C^{n \times n}$, its rank is defined as
\begin{align*}
\text{rank}(M)&=\text{number of linear independent columns of } M\\
&=\text{number of linear independent rows of } M.
\end{align*}
It should be noted that, for any integer $k<n$ the set
$
\{M \in C^{n \times n}: \text{rank}(M) \leq k \}
$
\textit{is not} a convex set.

\end{tcolorbox}

An important observation is that, by introducing  matrix $\WW$, the only source of nonconvexity is captured by the rank constraint. Indeed, as shown first in \cite{Lavaei:2012} and subsequently in \cite{madani2015convex},  in most cases this constraint can be dropped without affecting the OPF solution. 

To formally define the \textit{convexified version} of the AC-OPF problem, we note that bus voltage $\mathbf{V}$ appears in \emph{Voltage Limits} \eqref{eq:voltage constraints}
and \emph{Power Flow  equations} \eqref{eq:balance equations}. Therefore, these constraints are  redefined in terms of the new variable $\mathbf{W}$ as
\begin{subequations}\label{eq:voltage constraints W}
	\begin{align}
	\label{eq:voltage bound W}
	& (V_{k\,\min})^2\leq W_{kk}\leq (V_{k\,\max})^2, \quad \forall k \in \mathcal{N}\\
	\label{eq:voltage bound 2W}
	& W_{ll}+W_{mm}-W_{lm}-W_{ml}\le (\Delta V_{lm}^{\max})^2,\quad\forall(l,m)\in\mathcal{L}
	\end{align}
\end{subequations}
and 
\begin{subequations}\label{eq:balance equations W}
	\begin{align}
	\label{eq:active balance equation W}
	P^G_{k} -  P^L_{k}  =\sum_{\ell\in\mathcal{N}_k} \real\left\{(W_{kk}-W_{k\ell})^*Y_{k\ell}^*\right\}&, \quad \forall k \in \mathcal{N}\\
	\label{eq: reactive balance equation W}
	Q^G_{k} - Q^L_{k} =
	\sum_{\ell\in\mathcal{N}_k} \imag\left\{(W_{kk}-W_{k\ell})^*Y_{k\ell}^*\right\}&,\quad \forall k \in \mathcal{N}
	\end{align}
\end{subequations}
respectively.

\section{Control and State Variables}
\label{sec:control-state}

The  variables appearing in \eqref{eq:OPF Original} are usually divided into  two classes, depending on their role in the optimization problem. Indeed, already in \cite{Carpentier:1979aa}, see also the recent survey \cite{Capitanescu:2016}, the distinction between  control and state variables is explicitly made. As it is clear from their name, \emph{control variables} are those used by the network operator to set the operating condition of the network or, in other words, to \emph{control} its behavior. On the other hand, \emph{state variables} are \emph{dependent} variables that represent the state of a power network, and their values are a consequence of the designed control variables and of the constraints imposed to the power network. 

In particular, control and state variables are defined differently depending on the type of bus. In a generator bus (PV bus) $k\in\mathcal{G}$, (see e.g.\ \cite[Remark 1]{Low:2014a}) active power $P^G_{k}$ of the generator and magnitude $|V_k|$ of the complex bus voltage represent the control variables, 
while phase angle $\theta_k$ of bus voltage and generator reactive power $Q^G_{k}$ are the state variables. In a load bus (PQ bus), the active and reactive power of the load $P^L_{k},Q^L_{k}$ are given (their values are known to the network operator)
while magnitude and phase angle of bus voltage $|V_k|,\theta_k$ are state variables. A node to which both a generator and loads are connected is to be considered as a generator bus. 
Finally, in the slack bus $0$,  there is no control variable, while active and reactive  generator power $P^G_k$ and $Q^G_k$ are state variables.


To emphasize the inherent difference between control and state variables, we introduce the notation
\begin{align}\label{eq:control variables}
\uu&\doteq\{
P^G_1,  \ldots, P^G_{n_g},
|V_1|, \ldots, |V_{n_g}|\}\\
\label{eq:state variables}
\xx&\doteq \{Q^G_1,  \ldots ,Q^G_{n_g},
|V_{n_g+1}|,\ldots,|V_n|,\theta_1,\ldots,\theta_n\}
\end{align}
to denote respectively the control variables $\uu$ and the state variables $\xx$.
This allows to reformulate the problem \eqref{eq:OPF Original} in the following way

\noindent
\textbf{Reformulation of the Nominal AC-OPF}
\begin{eqnarray} \label{eq:OPF Restated}
&\underset{\uu}{\text{minimize}}
&f(\uu)\\ \nonumber
&\text{subject to:}& \text{there exist } \xx \text{ such that } g(\xx,\uu)=0 \text{ and }
h(\xx,\uu)\le 0, \nonumber 
\end{eqnarray}
where $f\doteq \sum_{k\in\mathcal{G}} f_k(P^G_{k})$,  the equality constraint $g(\xx,\uu)=0$ defines the Power Flow  Equations \eqref{eq:balance equations}, and the inequality constraint $h(\xx,\uu)\le 0$ summarize the power and voltage constraints \eqref{eq:power constraints} and \eqref{eq:voltage constraints}.

The formulation presented in \eqref{eq:OPF Restated} has the following interpretation: given the loads in the load buses, optimally design $\uu$ (active power of generators and voltage magnitude of generator nodes) such that \emph{there exists} a network state $\xx$ (reactive power and phase voltage at generator nodes $1,\ldots,n_g$, and complex voltage at load buses $n_{g+1},\ldots,n$) satisfying the operational constraints presented in \eqref{eq:balance equations}, \eqref{eq:power constraints}, and~\eqref{eq:voltage constraints}.

\section{Uncertain OPF}

The literature on OPF problems in the presence of load variations and uncertainties due to renewable energy resources (RES) has been constantly growing in the past years, testifying for an increasing interest in the problem. From a formulation viewpoint, the assumptions  on the network made in the classical formulation \eqref{eq:OPF Original} problem, i.e.\ (i)  completely predictable demand (and thus known in the design phase), (ii) no uncertainty in the amount of power being generated at  generator nodes, do not hold anymore  when dealing with modern power networks with high penetration of renewable sources and variable demands.

Hence, one needs to introduce a new category of generators, the \textit{renewable generators}, belonging to the set $\mathcal{R}\subseteq\mathcal{N}$ with $|\mathcal{R}|=n_r$\footnote{Hence, $\mathcal{N}=\{0\}\cup\mathcal{G}\cup\mathcal{R}\cup\mathcal{D}$}. 
A renewable energy generator connected to  bus $k\in\mathcal{N}$ provides an \textit{uncertain} complex power
\begin{equation}\label{eq: renewable uncertainty}
P^R_{k}(\delta^R_{k})+Q^R_{k}(\delta^R_{k})\i=P^{R,0}_{k}+Q^{R,0}_{k}\i+\delta^R_{k},
\end{equation}
with $P^{R,0}_{k}+Q^{R,0}_{k}\i$ being the nominal (predicted) power generated by the renewable energy source, and $\delta^R_{k}\in\boldsymbol{\Delta}^R_{k}\subset\mathbb{C}$ representing an uncertain complex fluctuation, which mainly depends on the environmental conditions, such as wind speed in the case of wind generators.  
Various works addressed the problem of properly modeling the power variation due to the presence of wind generators \cite{hodge_wind_2012,Hodge_forecasting}, and variable loads \cite{hodge2013short,hodge2012comparison}.

Similarly, the uncertain demand in bus $k\in\mathcal{N}$ is represented as
\begin{equation}\label{eq: load uncertainty }
P^L_{k}(\delta^L_{k})+Q^L_{k}(\delta^L_{k})\i=P^{L,0}_{k}+Q^{L,0}_{k}\i+\delta^L_{k}
\end{equation}
where $P^{L,0}_{k}$ and $Q^{L,0}_{k}$ denote the expected active and reactive load and  $\delta^L_{k}\in\boldsymbol{\Delta}^L_{k}\subset\mathbb{C}$ is the complex fluctuation in the demand at bus $k\in\mathcal{N}$.
The support set is the point $\{0\}$ if no uncertainty (i.e. no renewable generator or variable load) is present in bus $k$.

To simplify the notation, one may collect the different sources of uncertainty by introducing the \textit{uncertainty vector} 
\[
\dd \doteq
[\delta^L_{1}\,\cdots\,\delta^L_{n}\,\delta^R_{1}\,\cdots\,\delta^R_{n}]^T,
\]
which varies in the set 
\[
\DD \doteq \boldsymbol{\Delta}^L_{1} \times \cdots \times\boldsymbol{\Delta}^L_{n} \times
\boldsymbol{\Delta}^R_{1} \times \cdots \times \boldsymbol{\Delta}^R_{n}. 
\]

\subsection{Formulation of  robust and chance-constrained  OPF}

Once the sources of uncertainty have been modeled, the main idea in most of the literature is to propose solutions allowing to distribute the  power  mismatch among classical generator by introducing a \textit{deployment vector.} 
During real-time operation of the power network,  the amount of power mismatch---the difference between real-time (actual) and predicted demand---needs to be distributed among generators. This operation is called frequency control or primary and secondary control. In classical OPF, this is done through some coefficients which are generator specific. However, these coefficients are in general decided \textit{a-priori} in an ad-hoc fashion.

\begin{tcolorbox}[title= \sf \scriptsize \textbf{Deployment Vector},colframe=carmine!10,
colback=carmine!10,
coltitle=black,  
]
\sf\scriptsize

The  deployment vector physically represents the
way the frequency control is affecting the generator power.
It  was  seemingly  introduced in 
\cite{ETH-PMAPS} and \cite{ETH-Springer}, where it is referred to as ``distribution vector".
The same concept is present in many other works under different terminologies, such as  
 ``participation factor" in \cite{Jabr2013}, ``corrective control" in \cite{Jabr2015}, or ``affine control" in 
  \cite{Morari2013}.  The idea of   ``affinely adjustable robust counterpart (AARC)" adopted e.g.\ in \cite{Jabr2015}  can be also interpreted in a similar way.
 In many works,   the deployment vector was assumed to be constant, while more recent work explicitly consider  the possibility of optimizing over it so as to achieve a better
allocation of the generation resources (reserves). More sophisticated formulations have been subsequently introduced. For instance,  in \cite{ETH-TPS} a more
generic representation in the view of asymmetric reserves
is considered, while 
in \cite{ETH-Hawaii} the generation-load mismatch is redistributed also to some of the  loads.

\end{tcolorbox}

This approach worked well in cases where the amount of power mismatch  was not significant; however, once renewable generators are in the power network, this difference may become large, thus leading to line overloads in the network.
Modern approaches, as e.g.\  \cite{Bienstock2014,ETH-AC}, specifically incorporate these distribution parameters in the OPF optimization problem. To describe this idea, we formally introduce a 
\textit{deployment vector}
\[
\aa\doteq[\alpha_1,\ldots,\alpha_{n_g}]^T,
\]
with 
$\sum_{k\in\mathcal{G}}\alpha_k=1$, $\alpha_k \geq 0$ for all $k \in \mathcal{G}$,
whose purpose is to distribute among the available generators 
the power mismatch created by the uncertain generators and loads.

During operation, the active generation output of each generator is modified according to the  realization of the uncertain loads and RES power (which are assumed to be measured on-line) as follows
\begin{align}\label{eq: affine control}
\bar{P}^G_{k} &= P^G_{k}+\alpha_k\left(\sum_{j\in\mathcal{N}}\real\{\delta^L_{j}\}-\sum_{k\in\mathcal{R}}\real\{\delta^R_{k}\}\right)\\
&= P^G_{k}+\alpha_k \mathbf{s}^T \real\{\dd\},\quad \forall k \in \mathcal{G} \nonumber
\end{align}
with $\mathbf{s}^T\doteq[\mathbf{1}_n^T,\,-\mathbf{1}_n^T]$.
It is important to observe that, with the introduction of the uncertain renewable energy generator and load into the power network, summarized by the vector $\dd$, both Power Flow  equations \eqref{eq:balance equations} and Power Generation constraint \eqref{eq:power constraints} become uncertain.
In particular, the  equality constraints are rewritten as
\begin{subequations}\label{eq: uncertain balance equations}
	\begin{align}
	\label{eq: uncertain active balance equation}
	&P^G_{k}+\alpha_k \mathbf{s}^T \real\{\dd\}  + P^R_{k}(\dd) -  P^L_{k}(\dd) =
	\sum_{l\in\mathcal{N}_k} \real\left\{V_k(V_k-V_l)^*y_{kl}^*\right\}, \quad \forall k \in \mathcal{N}\\
	\label{eq: uncertain reactive balance equation}
	&Q^G_{k}  +  Q^R_{k}(\dd) - Q^L_{k}(\dd) =   \sum_{l\in\mathcal{N}_k} \imag\left\{V_k(V_k-V_l)^*y_{kl}^*\right\},\quad \forall k \in \mathcal{N},\\
	&\mathbf{1}^T\aa=1 \label{eq: uncertain balance equations alpha}
	\end{align}
\end{subequations}
and the power inequality constraints become
\begin{subequations}\label{eq: uncertain power constraints}
	\begin{align}
	\label{eq: uncertain power constraints active generation bound}
	&P_{k\,\min}\leq P^G_{k}+\alpha_k \mathbf{s}^T \real\{\dd\} \leq P_{k\,\max},\quad \forall k \in \mathcal{G} \\
	\label{eq: uncertain power constraints reactive generation bound}
	&Q_{k\,\min}\leq Q^G_{k} \leq Q_{k\,\max},\quad \forall k \in \mathcal{G},\\
	&\alpha_k \geq 0, \quad \forall k \in \mathcal{G}.
	\end{align}
\end{subequations}

Note that, following the classification introduced in Section \ref{sec:control-state}, the deployment vector $\aa$ can also be considered as a control variable, to be optimized to enhance performance. Hence, the set of control variables is redefined as
\begin{equation}\label{eq: control variables redefined}
\uu \doteq \{P^G_1,  \ldots, P^G_{n_g},
|V_1|, \ldots, |V_{n_g}|, \alpha_1,\ldots,\alpha_{n_g}\}.
\end{equation}
More importantly, it should be remarked that, with the introduction of uncertainties, the value of the state variables will depend on the specific realization of the uncertainty. Indeed, for every value of $\dd$, a different configuration of state variables may be necessary to satisfy the balance equation and other constraints. To emphasize this dependence on the uncertainty vector $\dd$, we use the following notation
 \begin{equation}
\xx_{\dd} \doteq \{Q^G_1,  \ldots, Q^G_{n_g},
|V_{n_g+1}|,\ldots,|V_n|,\theta_1,\ldots,\theta_n\}.
\end{equation}

This discussion allows to provide a formal statement of the  \emph{robust} version of the optimal power flow problem: denote by $g(\uu,\xx_{\dd},\dd) = 0$
the uncertain equality constraints collected in \eqref{eq: uncertain balance equations}, and by $h(\uu,\xx_{\dd},\dd) \le 0$ the uncertain inequalities collected in \eqref{eq: uncertain power constraints}.\\

\noindent{\bf Robust AC-OPF}
\begin{eqnarray} 
&\underset{\uu}{\text{minimize}} \
&f(\uu)
\label{eq: Robust AC-OPF}\\
\nonumber
& \text{subject to:}&
\text{for all } \dd \in\DD ,
\text{ there exist } \xx_{\dd} \text{ such that }\\ \nonumber
&&\qquad g(\uu,\xx_{\dd},\dd) = 0  \text{ and } h(\uu,\xx_{\dd},\dd) \le 0.
\end{eqnarray}

In this formulation of the robust OPF problem, the objective is to optimize the values of the nominal generated power, the bus voltage magnitude at generator node and the deployment vector so that i) the network operates safely for all values of the uncertainty and ii) the generation cost of the network is minimized.
In fact, if a solution to the problem above exists, we guarantee that for any admissible uncertainty, there exists a network state $\xx_{\dd}$ satisfying the operational constraints. 

In general, there are two paradigms to tackle uncertainty in optimization problems. The first approach is a deterministic worst-case approach where the constraints are enforced to remain feasible for the entire set of uncertainty. This approach is often not tractable except for rare cases where uncertainty enters the optimization problem in a ``simple'' fashion, e.g. affine, multi-affine, convex, etc. Furthermore, one can argue that such robust policies might be very conservative since some uncertainty scenarios are very unlikely to happen. 
The second approach is a probabilistic one where uncertainty is considered to be a random variable and constraints are enforced to hold for the entire set of uncertainty except a subset having arbitrary small probability measure.

In the deterministic worst-case approach of~\eqref{eq: Robust AC-OPF},  the constraints are enforced to hold for ``all'' possible values of the uncertain parameters. This  is in many cases excessive, and leads to  conservative results, with consequent degradation of the cost function (i.e., higher generation cost). 

Hence, to reduce conservatism, a chance-constrained approach is frequently adopted for OPF problems, in which a probabilistic description of the uncertainty is assumed to be  known, and a solution is sought which is valid for the entire set of uncertainty except for a (small) subset having probability  smaller than  a desired (small) risk level~$\varepsilon$. This approach is suitable for problems where ``occasional'' violation of constraints can be tolerated. One can argue that this is the case in power networks, since violation of line flow constraints does not necessarily lead to immediate line tripping. Rather, the line gradually heats up until a  critical condition is reached and only then the line is disconnected. Therefore, if line overload happens with low probability, this will not lead to line tripping nor it will damage the network. 

Formally, given a (small) risk level~$\varepsilon\in(0,1)$, the chance constrained version of the optimal power flow problem can be stated as follows\\

\noindent
\textbf{Chance-constrained AC-OPF}
\begin{eqnarray} 
&\underset{\uu}{\text{minimize}} \label{eq: OPF CCP}
&f(\uu)\\ \nonumber
& \text{subject to:}& 
\Pr\bigg\{\text{there exists a } \dd \in\DD\text{ for which does not exist }
 \xx_{\dd} \\ & & \ \ \ \ \ \ \ \ \ \ \ \ \text{  such that } 
 g(\uu,\xx_{\dd},\dd ) = 0  \text{ and } h(\uu,\xx_{\dd},\dd  ) \le 0 \bigg\}\leq \varepsilon. \nonumber
\end{eqnarray}	

Note that, in the above formulation, one accepts a certain \textit{risk} that the  designed control variable $\uu$ is such that some of the constraints might be violated for some value of the uncertainty, but this probability of violation is bounded by the \textit{a-priori} chosen violation level $\varepsilon$. Hence, one has a direct control on the risk of violating the constraints.

It should be remarked that the Robust AC-OPF problem in \eqref{eq: Robust AC-OPF}, and its probabilistic counterpart \eqref{eq: OPF CCP} 
are both computationally extremely hard, because the already nonlinear/nonconvex problem is made even more difficult by the
introduction of the semi-infinite robust constraints.   The presence of probabilistic constraints in  \eqref{eq: OPF CCP}  does not simplify the problem, since it requires the solution of hard multi-dimensional integration problems.  

\subsection{Solution approaches to  robust and chance-constrained  OPF}

The Robust AC-OPF formulation in~\eqref{eq: Robust AC-OPF}, with the distinction between independent (control) and dependent (state) variables,  was seemingly formulated in \cite{Zhang2011}, where an approach based on chance constrained programming is presented. The proposed solution is based on a  back-mapping approach and a linear approximation of the nonlinear model equations. The chance constrained counterpart is considered by~\cite{Wada:2014}, in which a non-convex iterative randomized method is provided, which partly mitigates the computational issues.

\subsubsection{DC-based approaches}

Most literature on uncertain OPF gets around the nonlinearities by recurring  to  the DC-based approximation discussed in Section~\ref{sec:DC}. These assumptions reduce the optimization problem to a quadratic program subject to uncertain linear equalities and inequalities, which still represents a challenging problems, at least for general probability distributions. 
First approaches in this direction have been based on scenario-tree generation methods, see for instance \cite{Yong2000}. These techniques suffer from severe computational complexity limitations, and do not offer theoretical guarantees on the probability of satisfaction of the constraints of the found solution.

In the case that uncertainties are assumed  to be Gaussian, the problem can be written in closed form \cite{Andersson2013}, or is amenable to a second-order cone-program \cite{Bienstock2014},  for which efficient solutions exist. This approach has been extended to ambiguous densities introducing the so-called robustified chance constraints in \cite{Lubin2016}. Similar ideas are at the basis of the approaches in \cite{Morari2015,Xie2018},  which consider distributionally robust approaches (that is, the 
solution has to remain valid for all uncertainties whose probability density  functions (pdfs) belong to a 
family of distribution functions sharing the same mean and variance, leading to the optimization of the so-called conditional value at risk  (CVar), which is again a convex problem.

The problem becomes much more difficult for general uncertainty distributions (indeed, it is known that the distribution of wind power is not Gaussian \cite{Hodge_forecasting}). In this case, a very promising approach is the application of recent results based on the so-called scenario approach, which are based on random  generation  of uncertainty samples. This is the approach followed for instance in \cite{ETH-PMAPS,ETH-TPS,ETH-Springer,
ETH-PSCC}.
In particular, in the field of the scenario based DC-OPF, a very promising approach has been recently proposed in~\cite{Modarresi2018}. 
In~\cite{Modarresi2018},  two approaches to solve DC-OPF via scenario are presented. First neglecting transmission constraints, it was shown that the number of i.i.d. samples does not grow with the size of the system. This makes the midnight solution of networks---when congestion does not exist but uncertainty from renewables is at the highest---much easier. This solution is based on a classical \emph{a-priori} scenario bound.
Then, by applying \emph{a-posteriori} scenario bounds, it was shown that one can start from a bad $\epsilon$ and progressively move to a good value of $\epsilon$ with a limited number of samples.
We refer the reader to the box ``Randomized Algorithms for Robust Optimization" for a more detailed discussion and pointers to the literature.

\subsubsection{Second-order and convex relaxations}

The limitations inherent to the DC approximation have motivated a few recent approaches to uncertain OPF, which employ more sophisticated relaxations able to capture also the reactive components of the power equations.
The work \cite{Perninge2013} makes use of second-order approximations of the stability boundary to approximate the probability of line violations. 
In \cite{ETH-AC} the SDP based relaxation of Section \ref{sec:relax} is exploited to derive a solution based on the scenario approach.
However, in order to guarantee solvability of the robust problem, the authors need to parameterize the dependence of the state variables on the uncertainty.
In that work, the authors cope with the need of guaranteeing the existence of a different value of $\WW$ for different values of the uncertainty $\dd$, by imposing a specific dependence on $\WW$ from the uncertainty. In our notation, \cite{ETH-AC} introduces the following finite (linear) parameterization
\begin{equation}
\label{eq:W-param}
\WW(\dd)=A+\sum_{k=1}^{n_r} B_k\dd_k
\end{equation}
where $A, B_1,\ldots,B_{n_r}$ become design variables in the optimization problem. 

\section{A randomized approach to chance constrained OPF}
\label{sec:scenario}

The approaches described in the previous section are surely conservative, since they either use a conservative DC-based approximation and/or   impose a very specific structure on the state variables. This section provides a possible way to overcome these limitations that uses recent results on randomized approaches to chance constrained optimization problems.

We start by noticing that in the convex relaxation of the AC-OPF problem discussed in Section~\ref{sec:relax}, some elements of $\WW$ involve the control variables $|V_{1}|,\ldots,|V_{n_g}|$ (specifically, the first $n_g$ entries in the diagonal of $\WW$), 
while  others are dependent variables corresponding to the voltage magnitude $|V_{n_g+1}|,\ldots,|V_n|$ at non-generator nodes, and the voltage phases $\theta_1,\ldots,\theta_n$. 
In order to distinguish between control and state variables appearing in it, \cite{OurTCNS2018} decomposes $\WW$ into the sum of  the two submatrices 
$\mathbf{W}^\mathbf{u} \doteq \text{diag}(|V_1|^2, \ldots,|V_{n_g}|^2, 0,\ldots),$
and
$\mathbf{W}^\mathbf{x} \doteq \mathbf{W}-\mathbf{W}^\mathbf{u}.
$
In this decomposition, we have a matrix $\WW^\mathbf{u}$ that includes the diagonal elements of $\WW$ corresponding to the generator nodes only, while the remaining elements of $\WW$ are collected in $\WW^\mathbf{x}$. With this in mind, the control and state variables are redefined as 
$
\uu \doteq \{\mathbf{P}^{G}, \aa,\WW^\uu\}
$
and 
$\xx_{\dd} \doteq \{\mathbf{Q}^G,\WW^\xx \},
$
respectively. 
With this notation settled, the convexified version of the robust AC optimal power flow problem is formally defined as follows \cite{OurTCNS2018}.

\noindent{\bf Convexified Robust AC-OPF (CR-AC-OPF)}
\begin{align}  \label{eq: CR-AC-OPF}
&\underset{\mathbf{P}^{G}, \aa,\WW_\uu}{\text{minimize }} 
\sum_{k\in\mathcal{G}} f_k(P^G_{k})\\ \nonumber
&\text{subject to: for all } \dd \in\DD ,
\text{ there exist } \mathbf{Q}^G=\mathbf{Q}^G(\dd),
\WW^\xx=\WW^\xx(\dd) \text{ such that } \\ \nonumber
&\quad \nonumber
\WW=\WW^\uu+\WW^\xx,
\quad  \nonumber
\mathbf{1}^T\aa=1,  \\
&\quad  \nonumber
P^G_{k}+\alpha_k \mathbf{s}^T \real\{\dd\} + P^R_{k}(\dd) -  P^L_{k}(\dd) =
\sum_{l\in\mathcal{N}_k} \real\left\{(W_{kk}-W_{kl})^*y_{kl}^*\right\},  \; \forall k \in \mathcal{N}\\
&\quad  \nonumber
Q^G_{k} +  Q^R_{k}(\dd) - Q^L_{k}(\dd) =
\sum_{l\in\mathcal{N}_k} \imag\left\{(W_{kk}-W_{kl})^*y_{kl}^*\right\}, \quad \forall k \in \mathcal{N}\\
&\quad  \nonumber
P_{k\,\min}\leq P^G_{k}+\alpha_k \mathbf{s}^T \real\{\dd\} \leq P_{k\,\max}, \quad \forall k \in \mathcal{G} \\
&\quad  \nonumber
Q_{k\,\min}\leq Q^G_{k}\leq Q_{k\,\max}, \quad \forall k \in \mathcal{G}\\
&\quad  \nonumber
\left(V_{k\,\min}\right)^2\leq W_{kk}\leq \left(V_{k\,\max}\right)^2,  \quad \forall k \in \mathcal{N}\\
&\quad  \nonumber
W_{ll}+W_{mm}-W_{lm}-W_{ml}\le (\Delta V_{lm}^{\max})^2,\ \quad\forall(l,m)\in\mathcal{L} \\
&\quad  \nonumber
\alpha_k\geq0,  \quad \forall k\in\mathcal{G} \nonumber
\quad  \nonumber
\WW\succeq 0.
\end{align}

\noindent

The problem above is a so-called 
robust optimization problem with  certificates, in which the design variables are those "controllable" by the network manager, i.e.
$
\uu \doteq \{\mathbf{P}^{G}, \aa,\WW^\uu\},
$
while the certificates are the variables that depend on the uncertainty  which are  "adjusted" to guarantee constraint satisfaction, i.e.
$\xx_{\dd} \doteq \{\mathbf{Q}^G,\WW^\xx \}.
$ Such a problem is very complex and extremely hard to solve. 

If one is willing to accept a (small) well defined risk of violation of the network constraints, then  a practical solution of such a problem can be obtained by using the randomized approach described in~\cite{SWC-CDC,SWC-TAC}. More precisely, 
assuming that independent and identically distributed (i.i.d.) samples of the uncertainty are available, we can formulate a problem that involves only those specific values of the uncertainty. In this case, although we cannot assure that the solution obtained satisfies the constraints, the risk of failure can be bounded \textit{a-priori}. This risk is a function of the number of samples used and can be made as small as desired; see~\cite{SWC-CDC,SWC-TAC}.
An important feature of scenario approach is that the underlying distribution of the uncertainty \emph{does not} need to be known. In fact, the only requirement is that the  uncertainty  is time-invariant, see e.g.\ \cite{campi_exact_2008,tempo_randomized_2013} for detailed discussion. 
Indeed,  when the distribution is known, i.i.d. random samples of it may be used, in a practical implementation one can use past observed values of the uncertainty obtained from observing the system's behavior.\\

\begin{tcolorbox}[title= \sf \scriptsize \textbf{Randomized Algorithms for Robust Optimization},colframe=carmine!10,
colback=carmine!10,
coltitle=black,  
]
\sf\scriptsize

The use of randomized algorithms for probabilistic analysis and design  has a long history.
The main ingredient in such techniques  is to  generate  a number of independent and identically distributed (i.i.d.) samples from the uncertainty set and then examine the performance function for all these random samples. Such an approach to probability estimation goes back to the works by Markov  \cite{markov_certain_1884}, Chebychev \cite{tchebichef1874valeurs},  Hoeffding \cite{hoeffding_probability_1963} and Bernstein \cite{Bernstein_1946}. This was later extended to the estimation of  the ``worst case" performance of a given performance index over the uncertainty set; e.g., see  \cite{tempo_probabilistic_1996,alamo_randomized_2012,alamo_randomized_2015,calafiore_probabilistic_2007,yasuaki_polynomial-time_2007,fujisaki_guaranteed_2007,dabbene_randomized_2010}.
Of particular interest to the problem addressed in this chapter is the use of randomized approaches to robust optimization. 
The so-called scenario approach to robust optimization  is a non-sequential algorithm which  is also based on extracting random samples from the uncertainty set and then solving an optimization problem subject to finite number of constraints defined by the samples drawn. This procedure was first introduced in \cite{calafiore_uncertain_2004} and later extended in \cite{campi_exact_2008,calafiore_scenario_2006,CamGar:11,calafiore_random_2010,vidyasagar_randomized_2001,alamo_randomized_2009,chamanbaz_sample_2013,chamanbaz_statistical_2014}. 

The results mentioned above can only deal with the case where there is a ``complete separation'' between the optimization variables and the uncertainty. In the case where some of the variables depend on the uncertainty, as it is the case in power networks, one needs an approach that can deal with so-called certificates. This problem was addressed in \cite{SWC-CDC,SWC-TAC} where the so-called  Scenario with Certificates (SwC) approach was developed.
	Consider the \textit{robust optimization problem with certificates} which is a generalization of the problem considered in this chapter
	\begin{eqnarray*}
	\min_{\uu} && f(\uu) \label{eq:certificates_opt}\\
	\text{subject to} && \forall \dd\in\DD  \ \exists \xx=\xx(\dd)
	\text{ satisfying } g(\uu,\xx,\dd)\le 0. \nonumber
	\end{eqnarray*}
	where $f(\uu)$ is a convex function. It was shown that problem above can be approximated by introducing the following \textit{scenario optimization  with certificates} problem, based  on the extraction of $N$ random samples $\dd^{(1)},\ldots,\dd^{(N)}$ of the uncertainty, 
	\begin{eqnarray*}
	\uu_{\text{SwC}}=\arg\min_{\uu,\xx_{1},\ldots,\xx_{N}}&& f(\uu)  \label{eq:scenario_cert}\\
	\text{subject to}&&  g(\uu,\xx_{i},\dd^{(i)}) \leq 0, \  i=1,\ldots,N.\nonumber
	\end{eqnarray*}
	Note that 
	in SwC \textit{a new certificate variable $\xx_{i}$ is created for every sample $\dd^{(i)}$}. In this way, one implicitly constructs an "uncertainty dependent" certificate, without assuming any \textit{a-priori} explicit functional dependence on $\dd$. 
    Let the number of samples $N$ be the smallest integer satisfying 
	\begin{equation*}
	\label{N_scenario}
	    \delta\leq \sum_{i=0}^{n_u-1}{N \choose i} \varepsilon^i(1-\varepsilon)^{N-i}
	\end{equation*}
		then, with probability $1-\beta$, the risk of violating the constraints of the robust optimization problem is less than or equal to $\varepsilon$; see \cite{SWC-CDC,SWC-TAC} for details. We remark that one can compute the exact number of samples by numerically solving above inequality which results in a non-linear equation; however, a sub-optimal bound on the number of samples can be obtained as follows
	    \begin{equation}
 		\label{eq: N_SwC}
		N\ge N\ped{SwC} = \frac{\mathrm{e}}{\varepsilon(\mathrm{e}-1)}\left(\ln\frac{1}{\beta}+n_{u}-1\right).
		\end{equation}
		
\end{tcolorbox}

Inspired by such an approach, we now describe a randomized algorithm that solves a modified version of the \textbf{CR-AC-OPF} problem where a small risk of failure is ``tolerated.'' The algorithm is based on  the generation of samples of the uncertainty $\dd$, and is guaranteed to return, with high probability, a solution guaranteeing the desired probabilistic guarantees. The procedure proposed in \cite{OurTCNS2018} is reported next.\\

\noindent{\bf SwC-AC-OPF Design Procedure}
\begin{itemize}
	\item[i)] Given probabilistic levels $\varepsilon$, and $\beta$ compute $N\ped{SwC}$ according to~\cite{SWC-CDC,SWC-TAC} (see note$^5$ on randomized algorithms for robust optimization).
	\item[ii)] Generate $N\ge N\ped{SwC}$ sampled scenarios $\dd^{(1)},\ldots,\dd^{(N)}$.
	\item[iii)]
	Solve the convex optimization problem {\bf SwC-AC-OPF} stated in 
	\eqref{eq: SwC Original}, which returns the control variables $\mathbf{P}^{G}, \mathbf{W^u},  \aa$.
	\item[iv)] During operation, measure  uncertainty in generations and loads $\dd$, and accommodate the $k$-th controllable generator as follows
\[
		\bar{P}^G_{k} =P^G_{k}+\alpha_k \mathbf{s}^T \real\{\dd\},
		\quad
		|V_{k}|  =\sqrt{W_{kk}},~k \in \mathcal{G}. 
\]

\end{itemize}

\noindent
In particular, the procedure involves the solution of the following random sampled problem, which is a convex semidefinite program solvable in polynomial time.

{\bf SwC-AC-OPF}
\begin{align}  \label{eq: SwC Original}
&\underset{
\tiny \begin{array}{l}
	\tiny \mathbf{P}^{G}, \mathbf{W^u},  \aa,\\
	\tiny \mathbf{Q}^{G,[1]},\ldots,\mathbf{Q}^{G,[N]},\\
	\tiny \WW^{\xx,[1]},\ldots, \WW^{\xx,[N]}
	\end{array}}{\text{minimize}} 
\gamma\\
&\text{subject to } \text{for } i=1,\ldots,N\nonumber \\
&\nonumber \WW^{[i]}=\WW^\uu+\WW^{\xx,[i]}\\
&\nonumber
L_{lm}^{[i]} =
|(W^{[i]}_{ll}-W^{[i]}_{lm})^*y_{lm}^*|
+
|(W^{[i]}_{mm}-W^{[i]}_{ml})^*y_{lm}^*|\\
&\sum_{k\in\mathcal{G}} f_k(P^G_{k}\!)+ \gamma_b \!\!\sum_{k\in\mathcal{G}}\!\! Q^{G,[i]}_{k}+ \!\gamma_\ell  \!\!\!\!\sum_{(l,m)\in\mathcal{L}^{\text{prob}}}\!\!   \!\!\!\!\!\!L_{lm}^{[i]} \le \gamma
\nonumber\\
&\nonumber
P^G_{k}+\alpha_k \mathbf{s}^T \real\{\dd^{(i)}\} + P^R_{k}(\dd^{(i)}) -  P^L_{k}(\dd^{(i)}) = 
\sum_{l\in\mathcal{N}_k} \real\left\{(W_{kk}^{[i]}-W_{kl}^{[i]})^*y_{kl}^*\right\},  \quad \forall k \in \mathcal{N}\\
&\nonumber
Q^{G,[i]}_{k} +  Q^R_{k}(\dd^{(i)}) - Q^L_{k}(\dd^{(i)}) = 
\sum_{l\in\mathcal{N}_k} \imag\left\{(W_{kk}^{[i]}-W_{kl}^{[i]})^*y_{kl}^*\right\}, \quad \forall k \in \mathcal{N}\\
&\nonumber
P_{k\,\min}\leq P^G_{k}+\alpha_k \mathbf{s}^T \real\{\dd^{(i)}\} \leq P_{k\,\max}, \quad 
Q_{k\,\min}\leq Q^{G,[i]}_{k}\leq Q_{k\,\max}, \quad \forall k \in \mathcal{G}\\
&\nonumber 
\left(V_{k\,\min}\right)^2\leq W_{kk}^{[i]}\leq \left(V_{k\,\max}\right)^2,  \quad \forall k \in \mathcal{N}\\
&\nonumber
W_{ll}^{[i]}+W_{mm}^{[i]}-W_{lm}^{[i]}-W_{ml}^{[i]}\le (\Delta V_{lm}^{\max})^2,\ \quad\forall(l,m)\in\mathcal{L} \\
&\nonumber
\alpha_k\geq0,  \; \forall k\in\mathcal{G}, 
\WW^{[i]}\succeq 0.\nonumber
\end{align}

We remark that the results in~\cite{SWC-CDC,SWC-TAC} provides an \emph{a-priori} guarantee that the \textbf{SwC-AC-OPF} design procedure is such that the probabilistic constraints of the Chance-Constrained AC-OPF \eqref{eq: OPF CCP} are satisfied with very high confidence. Moreover, although the number of variables in \textbf{SwC-AC-OPF} is proportional to the number of samples~$N\ped{SwC}$ (which can be high if one requires an extremely small risk of failure) the structure of the resulting optimization problem is amenable to parallelization and, hence, it can be efficiently solved even for large networks. Furthermore, sequential techniques such as the ones introduced in \cite{Chamanbaz_TAC_2016,Chamanbaz_CDC2013,calafiore_repetitive} can be used to reduce the computational complexity associated with solving \eqref{eq: SwC Original}. For example, in \cite{Chamanbaz_TAC_2016} at each iteration of the sequential algorithm, a sample complexity  smaller than the scenario bound is selected,  an optimization problem similar to \eqref{eq: SwC Original} but with smaller number of samples is solved and the robustness of the solution is checked in  a validation test. If the solution passes the validation test,  algorithm is successfully terminated; otherwise, a larger sample complexity is selected and a more complex optimization problem is solved; the two steps are repeated till a solution passes the validation test. Another approach with \emph{a-posteriori} probabilistic guarantees was recently introduced in \cite{Campi2018}, in which a desired sample complexity is chosen, the sampled-optimization problem is solved and then the  probabilistic robustness of the obtained solution is evaluated by computing the number of support constraints.   

\section{Conclusion}
In this chapter, we reviewed a number of methodologies for solving the problem of Optimal Power Flow (OPF) in the presence of renewable energy sources (RES) and uncertain loads. The presence of RES and uncertain loads is inevitable in the modern power grids. Their presence introduces huge amount of uncertainty into the grid which needs to be taken into account when dispatching the controllable (conventional) generators, otherwise, network constraints are very likely to be violated leading to line tripping or even cascading failures. 
In solving this problem one not only needs to cope with uncertainty but also needs to address non-convexity associated with the presence of power flow equations which are quadratic equality constraints. 
For this reason, we first reviewed state of the art  relaxation techniques to circumvent non-convexity in the formulation and next reviewed robust techniques applicable to this problem in order to handle uncertainty. In particular, we focused our attention to some randomized methods recently developed capable of efficiently solving the problem.

\bibliographystyle{plain}
\bibliography{ref_frg}
\end{document}